\documentclass[letterpaper, 10 pt, conference]{ieeeconf}  
\IEEEoverridecommandlockouts
\usepackage{balance}
\usepackage{graphicx}
\usepackage{algorithm2e}
\usepackage{algpseudocode}
\RestyleAlgo{ruled}
\usepackage{ifthen}                                                
\usepackage{cite}
\usepackage{amssymb,amsmath}

\renewcommand{\epsilon}{\varepsilon}
\renewcommand{\phi}{\varphi}
\renewcommand{\d}{\,\mathrm{d}}

\newcommand{\R}{\mathbb{R}} 

\title{\LARGE \bf
Optimization of External Stimuli for Populations of Theta Neurons\\ via Mean-Field Feedback Control*
}

\author{Roman Chertovskih$^{1}$, Nikolay Pogodaev$^{2}$, Maxim Staritsyn$^{1}$, Joaquim Da
Silva Sewane$^{3}$\\ and Ant\'onio Pedro Aguiar$^{1}$
\thanks{*The authors acknowledge the financial support of the Foundation for Science and
Technology (FCT, Portugal) 
in the framework of the Associated Laboratory ``Advanced Production and Intelligent
Systems'' (AL ARISE, ref. LA/P/0112/2020), R\&D Unit SYSTEC (base UIDB/00147/2020 and programmatic UIDP/00147/2020
funds), and projects SNAP (ref. NORTE-01-0145-FEDER-000085) and MLDLCOV (ref. DSAIPA/CS/0086/2020).
}
\thanks{$^{1}$Roman Chertovskih, Maxim Staritsyn and Ant\'onio Pedro Aguiar are with 
Research Center for Systems and Technologies (SYSTEC), 
Faculty of Engineering, 
University of Porto, 
Rua Dr. Roberto Frias, 
s/n 4200-465,
Porto, Portugal
{\tt\small roman@fe.up.pt, staritsyn@fe.up.pt, pedro.aguiar@fe.up.pt}
}%
\thanks{$^{2}$Nikolay Pogodaev is with 
Department of Mathematics ``Tullio Levi-Civita'', 
School of Sciences, 
University of Padova, 
Via Trieste, 63 - 35121 Padova, Italy
{\tt\small nickpogo@gmail.com}
        }%
\thanks{$^{3}$ Joaquim Da Silva Sewane is with 
Department of Mathematics and Informatics, 
Faculty of Sciences, 
University of Eduardo Mondlane, 
Av. Julius Nyerere, nr. 3453 Maputo, Mozambique
{\tt\small joaquimdasilvasewane@gmail.com}
}%
}
\begin{document}

\maketitle
\thispagestyle{empty}
\pagestyle{empty}

\begin{abstract}

We study a problem of designing ``robust'' external excitations for control and
synchronization of an assembly of homotypic harmonic oscillators representing so-called
theta neurons. The model of theta neurons (Theta model) captures, in main, the bursting
behavior of spiking cells in the brain of biological beings, enduring periodic
oscillations of the electric potential in their membrane. 

We study the following optimization problem: to design an external stimulus (control),
which steers all neurons of a given population to their desired phases (i.e.,
excites/slows down its spiking activity) with the highest probability. 

This task is formulated as an optimal mean-field control problem for the local continuity
equation in the space of probability measures. To solve this problem numerically, we
propose an indirect deterministic descent method based on an exact representation of the
increment (infinite-order variation) of the objective functional. We discuss some aspects
of practical realization of the proposed method, and provide results of numerical experiments.

\end{abstract}

\section{INTRODUCTION}

The phenomenon of synchronization of oscillatory processes arise
in many physical and natural systems involving (relatively large) collections of structurally similar interacting objects. 
This type of behavior~--- typically manifested in practice
by a formation of (desired or pathological) time-periodic patterns~--- is demonstrated, e.g., by semiconductors in laser physics \cite{Fisher-Laser}, vibrating processes in mechanics \cite{blekhman1988synchronization}, biochemical reactions \cite{Nishikawa,kuramoto2003chemical}, as well as in cardiac and neural activity \cite{UHLHAAS2006155,Schiff1994,Glass}. 

In connection with oscillatory processes, there naturally arise problems of designing 
artificial signals that can drive open systems towards (or away from) synchronous oscillations and frequency entrainment; important examples 
are clinical treatment of neurological and cardiac deceases (such as Parkinson’s disease, epilepsy, and cardiac arrhythmias), control of circadian rhythms \cite{winfree2013geometry}, organization/destruction of patterns in complex dynamic structures \cite{Istvan}, and in neurocomputing \cite{Hoppensteadt2000,Hoppensteadt2001}. 

Starting from the pioneer works of Y. Kuramoto and H.~Araki, the mathematical imperative in the study of oscillatory ensembles is the \emph{mean field} dynamics, which describes the behavior of an ``averaged'' representative of the population instead of tracking all individuals in person. This approach leads to a treatable (and elegant) mathematical representation of the ensemble dynamics even in the case when the cardinality of the population becomes very large, and is naturally translated to the control-theoretical context: 
in the most of applications, 
it is technically difficult (or even impossible) to ``isolate'' the \emph{control influence}  for a particular oscillatory unit; on the contrary, admissible signals usually affect a significant part of the system, or the system as a whole. 
The topic of control engineering which is focused on designing ``simultaneous'' control signals for multi-agent systems is familiar under the name \emph{ensemble control}. ``Adaptive'' (distributed in the phase space) signals are called \emph{mean-field type} controls. 

In this paper, we address a particular optimal control problem of the type \cite{Li-TACON-2013} based on a classical oscillatory model \cite{Theta} from the mathematical neuroscience. Namely, we study the problem of in-phase synchronization of the mean field of so-called theta neurons: to steer a given probability distribution of harmonic phases 
towards a target one by a simultaneous (ensemble) or individual (mean-field) control. 

To solve our problem numerically, we propose a deterministic iterative method of sequential ``control improvement'', entailed by an an exact formula for the variation of the objective functional.  The proposed approach is based on the optimal mean-field control theory (the dynamic optimization in the space of probability measures) and is quite flexible: it admits one to treat arbitrary statistical ensembles, and can be applied to any problem of a ``state-linear'' structure, far beyond the considered specific model.

\section{PROBLEM STATEMENT.\\ MEAN-FIELD CONTROL SETUP}

Consider a population of homotypic oscillatory systems represented by the canonical Ermentrout-Kopell model  \cite{Theta,Strogatz}. This model 
describes the time-evolution of excitable neurons (customary named ``theta neurons'') which endure periodic oscillations of their membrane potential. Each theta neuron in the population is characterized by its phase 
\[
\theta(t) \in \mathbb S^1\doteq \R/2\pi\mathbb Z 
\]
which satisfies the ODEs
\begin{equation*}
    \label{Theta}
    \frac{\d }{\d t}\theta \doteq \dot \theta = v_{u}(\theta, \eta) \doteq (1-\cos\theta) +(1+\cos\theta)\left(u+\eta\right).
\end{equation*}
Here, $\eta$ is the baseline current in the neuron membrane, which varies in a given interval $\mathbb I \doteq [a, b]$, 
and $u$ is an external stimulus. 

Theta model provides a simple mathematical description of the so-called  \emph{spiking} behavior. By convention, we say that a neuron produces a spike 
at time $t$ if $\theta(t)=\pi$.  If $\eta>0$ (and $u \equiv 0$) the neuron spikes periodically with the frequency $2 \sqrt{\eta}$. If $\eta<0$, the neuron is excitable and can produce spikes after a sufficiently intensive stimulus $u$. 

In what follows, $\eta$ is viewed as a parameter of the model fluctuation. 
In the simplest case, this parameter runs through a finite set 
\(
\{\eta_k, \ k=\overline{1,N}\},
\) 
which corresponds to a finite ensemble 
\(
\{\theta_k, \ k=\overline{1,N}\}
\)
of theta neurons,
\begin{equation}
\dot \theta_k = v_{u}(\theta_k, \eta_k), \quad k=\overline{1,N}.\label{Theta-k}    
\end{equation}
In a more general setup to be discussed below, $\eta$ can be drawn from a given probability distribution.

Remark that \eqref{Theta-k} falls into the well-recognized Watanabe-Strogatz class of 
phase oscillators 
driven by complex functions $t \mapsto H_k(t) \in \mathbb C$, 
\[
\dot \theta_k = \omega_k + {\rm Im}\big(H_k(t)\, e^{-i \, \theta_k}\big), \quad k=\overline{1,N},
\]
where $\omega_k$ is the natural (intrinsic) frequency of the $k$th oscillator in the population, and  $H_k$ is the associated input, modulated by a sinusoidal function (sometimes, 
this 
model is called ``sinusoidally coupled''); in general, both the natural frequencies and
the inputs can be effected by an external driving parameter, furthermore, $H_k$ can model
interactions between oscillators inside the population. Note that model \eqref{Theta-k} fits the general statement with
\[
\omega_k=\omega_k(u) \doteq u + \eta_k + 1,\]
\[
H_k =H_k(u) \doteq i(u + \eta_k - 1),
\]
which does not involve interaction terms (formally, equations  \eqref{Theta-k} are paired
only by the common term $u$). In the context of applications, this non-interacting model
can be viewed as a ``first-order approximation'' of a sufficiently sparsely connected
neural network (such are real biological ones), especially, if the neurons' activity is
studied over relatively short time periods. The case of interacting
neurons will be briefly discussed in section \ref{sec-conc}.

\subsection{Mean-Field Limit}

We are interested in the behavior of system \eqref{Theta-k} for the case when $N\to \infty$. Introduce extra, ``fictitious'' states $t \mapsto \eta_k(t)$ as solutions to
\begin{equation}
\dot \eta_k =0, 
\end{equation}
accompanying \eqref{Theta-k}, and consider the empirical probability measure 
\begin{equation}
\label{mu-N}
\mu_t^N =\frac{1}{N} \sum_{k=1}^N \delta_{(\theta_k(t), \eta_k(t))}, 
\end{equation}
($\delta_x$ stands for the Dirac probability measure concentrated at at a point $x$). 

The measure-valued function $t \mapsto \mu_t^N$ designates the statistical behavior of the ensemble $\{(\theta_k,\eta_k), \ k=\overline{1,N}\}$: for any Borel set $A \subset \mathbb S^1\times\mathbb I$, the value $\mu^N_t(A)$ shows the number of neurons whose phase belongs to $A$. 

It is well-known that the curve $t \mapsto \mu_t^N$ satisfies, in the weak sense, the local continuity equation \cite{AmbrosioBook2007}
\begin{equation}\label{PDE}
\partial_t\mu_t(\theta, \eta) + \partial_\theta\big(v_{u}(\theta, \eta)\, \mu_t(\theta, \eta)\big) = 0. 
\end{equation}
Recall that the map $t \mapsto \mu_t$ is said to be a weak (distributional) solution of~\eqref{PDE} iff
\begin{equation*}
\label{conteq-sol}
\begin{array}{c}
0 = \displaystyle\int_0^T \d t \int_{\mathbb S^1\times\mathbb I} \big(\partial_t \phi + \nabla_x \phi \cdot v_u\big)\d \mu_t\\[0.4cm]
\forall \, \phi \in C^1_c((0,T)\times \mathbb S^1\times\mathbb I).
\end{array}
\end{equation*}
($C^1_c((0,T)\times \mathbb S^1\times\mathbb I)$ denotes the space of continuously differentiable functions $(0,T)\times \mathbb S^1\times\mathbb I \mapsto\R$ with compact support in $(0,T)\times \mathbb S^1\times\mathbb I$.) Under standard regularity assumptions, the weak solution exists, it is unique, and it is absolutely continuous as a function $[0, T]\mapsto \mathcal P(\mathbb S^1\times \mathbb I)$; here $\mathcal P(\mathbb S^1\times \mathbb I)$ denotes the space of probability measures on $\mathbb S^1\times \mathbb I$ endowed with any Wasserstein distance $W_p$, $p\geq 1$ \cite{AmbrosioBook2007}.

Equation \eqref{PDE} provides the macroscopic description of the population of microscopic dynamical units \eqref{Theta-k} called the \emph{mean field}. This representation remains valid in the limit $N\to \infty$, when $\mu^N$ converges to some $\mu \in \mathcal P(\mathbb S^1 \times \mathbb I)$ in $C([0, T];\mathcal P(\mathbb S^1 \times \mathbb I))$.
Moreover, \eqref{PDE} makes sense if phases $\theta$ and currents $\eta$  are drawn from an abstract probability distribution on the cylinder $\mathbb S^1 \times \mathbb I$,
\begin{equation}\label{PDE-0}
\mu_0=\vartheta \in \mathcal P(\mathbb S^1 \times \mathbb I).
\end{equation}
Indeed, one can immerse the system of ODEs \eqref{Theta-k} in a deterministic $(\mathbb S^1 \times \mathbb I)$-valued random process
\[
(t,\omega)  \mapsto \Theta_t(\omega),
\]
defined on a probability space $(\Omega, \mathcal F, \mathbb P)$ of an arbitrary nature ($\Omega$ is an abstract set, $\mathcal F$ is a sigma-algebra on $\Omega$, and $\mathbb P$ is a probability measure $\mathcal F\mapsto [0,1]$), and satisfying the ODE 
\[
\frac{\d }{\d t}\Theta_t(\omega) = 
\left(\begin{array}{c}
v_u\big(\Theta_t(\omega)\big)\\[0.2cm] 
0 
\end{array}\right).
\]
It is a simple technical exercise to check that the function 
\[
t\mapsto \mu_t \doteq (\Theta_t)_\sharp \mathbb P
\]
solves the Cauchy problem \eqref{PDE}, \eqref{PDE-0} with $\vartheta \doteq (\Theta_0)_\sharp \mathbb P$, where the symbol $\sharp$ denotes the operation of pushforward of a measure by a (Borel) function $\Omega \mapsto \mathbb S^1 \times \mathbb I$. Note that empirical ensembles \eqref{mu-N} fit this setup if $\Omega=\{1,\ldots, N\}$ and $\mathbb P$ is the normalized counting measure.

Finally, observe that the variable $\eta$ enters PDE \eqref{PDE} as a parameter rather
than state variable. This means that \eqref{PDE} can be regarded as an $\eta$-parametric
family of continuity equations on the 1D space $\mathbb S^1$ rather than a PDE on the 2D
space $\mathbb S^1\times\mathbb I$. This observation is essential for the numerical
treatment of the problem \eqref{PDE} (see section~\ref{sec-num}).

\subsection{Control Signals}

Now, we shall fix the class of admissible control signal $u$. Consider two options: 
\begin{itemize}
\item $u=u(t)$, i.e., the control effects all neurons of the ensemble in the same way. 
We call this type of external influences the \emph{ensemble} (simultaneous, common) \emph{control}. Such a control is statistical in its spirit as it influences the whole ensemble ``in average''. As a natural space of such controls we choose
\begin{equation}
 u \in \mathcal U \doteq L_2([0,T]; \R).\label{u}
\end{equation}

\item $u=w_t(\theta, \eta)$, i.e.,  
the stimulus is adopted to the neuron's  individual characteristics and phase-dependent. The use of such a distributed, \emph{mean-field type control}
\begin{equation}
 w \in {\mathcal W} \doteq L_2([0,T]; C(\mathbb S^1\times \mathbb I;\R)),\label{w}
\end{equation}
assumes some technical option to variate control signals over the spatial domain.  
\end{itemize}
It is natural to expect that the second-type control should perform better. However, let us stress again that the practical implementation of ``personalized'' control signals is hardly realistic as soon as the number of driven objects is large enough (for experiments that pretend to mimic the biological neural tissue, this number should be astronomic!). 
In reality, a meaningful class of control signals is $\mathcal U$, or something ``in the middle'' between the mentioned two options.

\subsection{Performance Criterion}

We study a generalization of the optimization problem \cite{Li-TACON-2013}: to steer the neural population to a target phase distribution at a prescribed (finite)  time moment $T>0$ with care about the total energy of the control action. Assuming that the target distribution is given by a (bounded continuous) function $\eta \mapsto \check \theta(\eta)$, our optimization problem reads:
\[
(P_1) \qquad 
\left\{
\begin{array}{c}
\min I[u] = \displaystyle 
\int F\big(\theta,\check\theta(\eta)\big)\d \mu_T(\theta, \eta)\\[0.3cm]
\ \ \qquad\displaystyle +\frac{\alpha}{2}\int_0^T u^2(t) \d t, \quad \alpha>0,\\[0.5cm]
\displaystyle\mbox{subject to \eqref{PDE}, \eqref{u}},
\end{array}
\right.
\]
where
\begin{align*}
F(\theta,\omega)= & \, \frac{1}{2}(\sin \theta-\sin \omega)^2+\frac{1}{2}(\cos \theta -\cos \omega)^2\\
= &
1 - \cos(\theta-\omega),
\end{align*}
and 
\[
\displaystyle\int \doteq \int_{\mathbb S^1\times\mathbb I}.
\] 
In this problem, the part of state variable is played by the probability measure $\mu_t$. 

Note that the functional $I$ and the dynamics \eqref{PDE} are linear in $\mu$ (despite the non-linearity of the map $(\theta, \eta) \mapsto v_u(\theta, \eta)$). At the same time, \eqref{PDE} contains a product of $\mu$ and $u$, which means that $(P_1)$ is, in fact, a \emph{bi-linear} (non-convex) problem. 

Standard arguments from the theory of transport equations in the Wasserstein space \cite{AmbrosioBook2007} together with the classical Weierstrass theorem ensure that problem $(P_1)$ is well posed, i.e., it does have a minimizer within the admissible class $\mathcal U$ of control signals (refer, e.g., to \cite{POGODAEV20203585}). 

An alternative version of problem $(P_1)$ is formulated in terms of the mean-field type control: 
\[
(P_2) \quad 
\left\{
\begin{array}{c}
\min J[w]= \displaystyle \int F\big(\theta,\check\theta(\eta)\big)\d \mu_T\\[0.3cm] \qquad \qquad + \displaystyle\frac{\alpha}{2}\int_0^T\!\!\! \d t \int w_t^2 \d \mu_t,\\[0.5cm]
\displaystyle\mbox{subject to \eqref{PDE}, \eqref{w}}.
\end{array}
\right.
\]
In what follows, we shall focus on the ``more realistic'' statement $(P_1)$, though all the forthcoming results can be extended, at least formally, to problem $(P_2)$.

\section{COST INCREMENT FORMULA.\\ NUMERICAL ALGORITHM}

As it was remarked above, problem $(P_1)$ is linear in state-measure. This fact allows us to represent the variation of the cost functional $I$ with respect to \emph{any} variation of control $u$ \emph{exactly} (without any residual terms). The announced representation follows from the duality with the co-state from Pontryagin's maximum principle \cite{Pogo-ContEq}, and generalizes the classical exact increment formula for conventional state-linear optimal control problems \cite{SChPP-2022}.

Consider two arbitrary controls 
\[
\bar u, u \in \mathcal U, \quad u \neq \bar u, 
\]
and let 
\[
t \mapsto \bar\mu_t\doteq \mu_t[\bar u]\mbox{ and }t \mapsto \mu_t\doteq \mu_t[u]
\]
be the respective weak solutions to the continuity equation \eqref{PDE}. Let also 
\[
\bar p\doteq p[\bar u]: \ (t, \theta, \eta) \mapsto \bar p_t(\theta, \eta)
\] 
be a classical solution to the following (non-conservative transport) equation:
\begin{align}
\partial_t p_t(\theta, \eta)  + & \ \partial_\theta p_t(\theta, \eta) \cdot v_{\bar u(t)}(\theta, \eta) = 0\label{p}.
\end{align}
PDE \eqref{p} is known to be dual to the (conservative transport equation) \eqref{PDE}; the duality is formally established by
the observation that 
the map \[t \mapsto \displaystyle\int \bar p_t \d \bar \mu_t\] is constant on $[0,T]$. 
One can check that, under the common regularity of the problem data, this map is an absolutely continuous function $[0,T] \mapsto \R$  (refer to \cite{AmbrosioBook2007} for further details).

As soon as $\bar p$ is chosen as a solution to \eqref{p} with the terminal condition
\begin{align}
p_T(\theta, \eta) = & -F\big(\theta,\check\theta(\eta)\big),
\label{p-T}
\end{align}
the discussed duality makes it possible to represent the increment (variation)
\[
\Delta I \doteq I[u] - I[\bar u]
\]
of the functional $I$ as follows:
\begin{align}
    - 
    \Delta I  = 
    \int_0^T \big(H\left(\mu_t, \partial_\theta \bar p_t, u(t)\right) -  H\left(\mu_t, \partial_\theta \bar p_t, \bar u(t)\right)\big)\d t,\label{incr}
\end{align}
where
\[
H(\mu, \zeta, u) \doteq u \int \zeta(\theta, \eta) \cdot (1 + \cos \theta) \d \mu(\theta, \eta) - \frac{\alpha}{2} u^2.
\]
The derivation of this formula is dropped, since it is completely similar to \cite{SChPP-2022}.
 
 Based on representation \eqref{incr}, we can treat problem $(P_1)$ in the following iterative way: given a reference control $\bar u$, one looks for a new ``target'' signal $u$ that ``improves'' the functional value, i.e such that $\Delta I < 0$. The best choice of the target control is provided by the maximization of the integrand of \eqref{incr} in the variable $u$:
\[
H\left(\mu_t, \partial_\theta \bar p_t, u\right) \to \max, \quad u \in \mathbb R.
\]
The unique solution of the latter problem is obtain in the analytic form as
\begin{equation}\label{u-mu}
u_t[\mu] = \frac{1}{\alpha}\int \partial_\theta \bar p_t(\theta, \eta)\, (1 + \cos \theta) \d \mu(\theta, \eta).
\end{equation}
Here, it is worthwhile to mention that the reference dual state $\bar p$ enters formula \eqref{u-mu} only in the form of the partial derivative 
\[
\bar \xi_t(\theta, \eta) \doteq \partial_\theta \bar p_t(\theta, \eta).
\]
Differentiating \eqref{p} and \eqref{p-T} in $\theta$ one can easily check that $\bar \xi$ solves the $\eta$-parametric family of the same continuity equations \eqref{PDE} \emph{backward in time}, starting from the terminal condition
\begin{equation}
\label{dual-term}
\xi_T =-\partial_\theta F\big(\theta,\check\theta(\eta)\big) \doteq \sin\left(\check\theta(\eta)-\theta\right).
\end{equation}
Now, \eqref{u-mu} can be reformulated in terms of the variable $\bar \xi$:
\begin{equation}
\label{u-1}
u_t[\mu] = \frac{1}{\alpha}\int \bar \xi_t(\theta, \eta)\, (1 + \cos \theta) \d \mu(\theta, \eta).
\end{equation}

Note that the map $(t, \mu) \mapsto u_t[\mu]$ can be used as a feedback control 
\[
[0, T]\times \mathcal P(\mathbb S^1\times\mathbb I) \mapsto \R 
\]
of system \eqref{PDE} in the space of probability measures. Injecting this control into \eqref{PDE}, we obtain a nonlocal continuity equation 
\begin{equation}
\partial_t\mu_t + \partial_\theta\big(v_{u[\mu_t]}\, \mu_t\big) = 0, \quad \mu_0=\vartheta,\label{PDE-back}
\end{equation}
which is well-posed (thanks to the fact that function $(\theta, \eta) \mapsto v_u(\theta, \eta)$ is smooth and bounded). 
Solving the last equation numerically, and substituting its solution $t \mapsto \hat\mu_t\doteq \hat \mu_t[\bar u]$ into \eqref{u-mu}, 
we construct the ``improved'' signal: 
\[
u(t) = u_t[\hat \mu_t].
\]

This idea gives rise to the following Algorithm~\ref{alg}. 
\begin{algorithm}
\caption{Numerical algorithm for optimal ensemble control}
\label{alg}
\KwData{$\bar u \in \mathcal U$ (initial guess), $\varepsilon>0$ (tolerance)}
\KwResult{$\{u^k\}_{k \geq 0} \subset \mathcal U$ such that $I[u^{k+1}] < I[u^{k}]$}
$k \gets 0$\;
$u^0 \gets \bar u$\;
\Repeat{$I[u^{k-1}] - I[u^{k}] < \varepsilon$}{
$\mu^{k} \gets \hat \mu[u^k]$\;
$u^{k+1} \gets u[\mu^k]$\;
$k \gets k+1$\;
  }
\end{algorithm}

By construction, Algorithm~\ref{alg} generates a sequence $\{u^k\}_{k \geq 0} \subset \mathcal U$ of controls with the property: 
\[
I^{k+1}\doteq I[u^{k+1}] < I[u^{k}]\doteq I^k.
\]
Since the sequence of numbers $(I^k)_{k\geq 0}$ is bounded from below by $\min(P)$ it converges.

Finally, remark that the same line of arguments can be formally applied to problem $(P_2)$. The respective mean-field type control takes the form
\[
w_t(\theta, \eta) = \frac{1}{\alpha}\bar \xi_t(\theta, \eta)\, (1 + \cos \theta).
\]
This construction gives rise to an iterative method, similar to Algorithm~\ref{alg}. 

\section{NUMERICAL RESULTS}\label{sec-num}

Let us discuss several aspects of the numerical implementation of Algorithm~\ref{alg}. 

First, note that the method proposed here does not involve any intrinsic parametric optimization: the
most of indirect algorithms for optimal control require the dynamic adjustment of some
internal computational parameters; such are standard methods based on Pontryagin's maximum
principle \cite{Arguchintsev2009,Drag} that imply the internal
such as line search for the specification of the ``depth'' of the needle-shaped (or weak)
control variations. 

Each iteration of Algorithm~\ref{alg} requires numerical solution of two problems: 
one is the linear problem \eqref{PDE}, \eqref{dual-term} (integrated backward
in time), and one for the nonlocal continuity equation \eqref{PDE-back} (solved
numerically forward in time). Since both \eqref{PDE} and \eqref{PDE-back} have no terms involving partial derivatives 
in $\eta$, one can think of $\eta$ as a parameter and solve the corresponding parametric 
families of one-dimensional continuity equations.

\medskip

Consider the problem $(P)$ with initial distribution of neurons $\mu_0$ given by the density function 
\[
\rho_0(\theta,\eta) = \big(2+3\cos(2\theta)-2\sin(2\theta)\big)\eta, 
\]
and with constant target function $\check\theta(\eta) \equiv \pi$.
In other words, our goal is to bring neurons' states as close as possible to the segment $0\times \mathbb I$ by the time moment $T$ with the aid of sufficiently small controls.

Parameters for the computation: 
\[
T=6, \quad \mathbb I=[0.0,1.0], \quad \alpha=1; 
\]
we used 512 Fourier harmonics in $\theta$ and grid steps   
\[
\Delta \eta = 0.002, \quad \Delta t = 0.002.
\]
Equations \eqref{PDE} and \eqref{PDE-back} are integrated by the standard spectral method \cite{Boyd2001} using the trigonometric Fourier expansion in $\theta$ for each $\eta$ from the grid. Parameters of the algorithm: $\bar u \equiv 0$, $\varepsilon=0.01$.

\begin{figure}[!ht]
  \centerline{\includegraphics[width=0.5\textwidth]{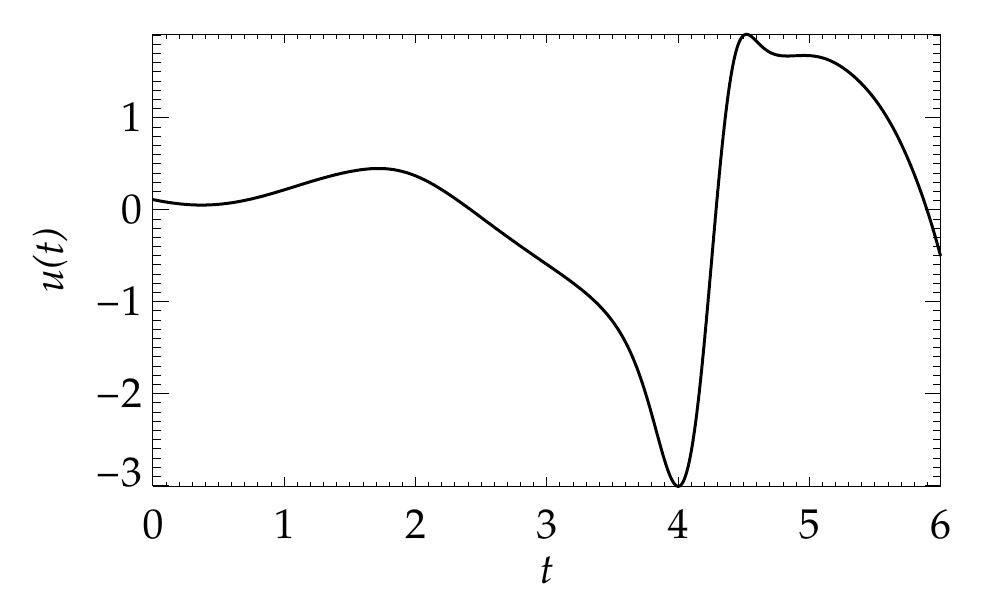}}
  \caption{Control input computed by the Algorithm~\ref{alg}}
  \label{fig:ensemble_control}
\end{figure}

\begin{figure}[!ht]
  \centering
  \includegraphics[width=0.40\textwidth]{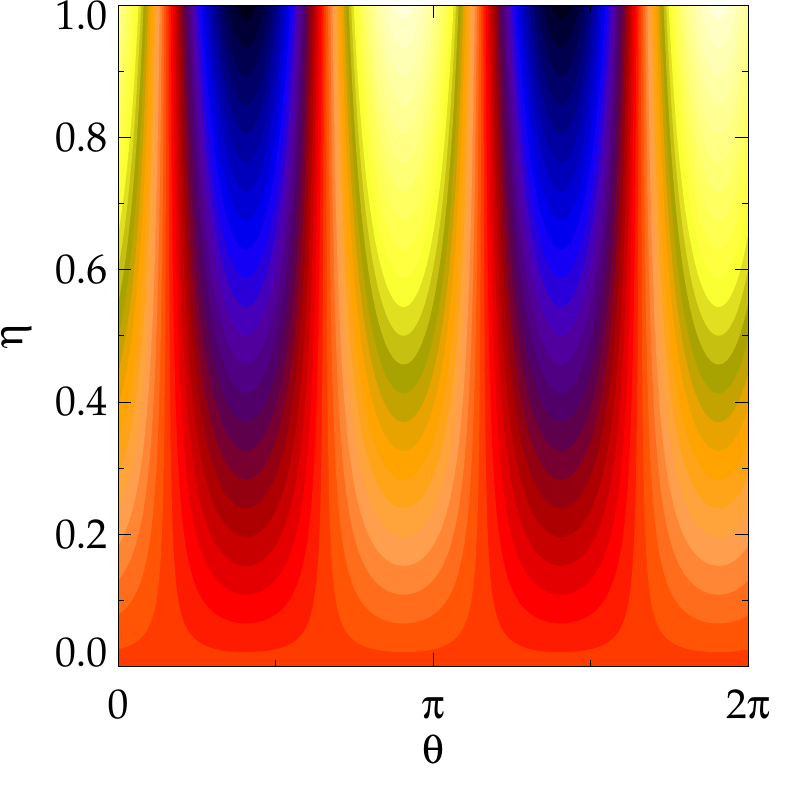}
  \includegraphics[width=0.40\textwidth]{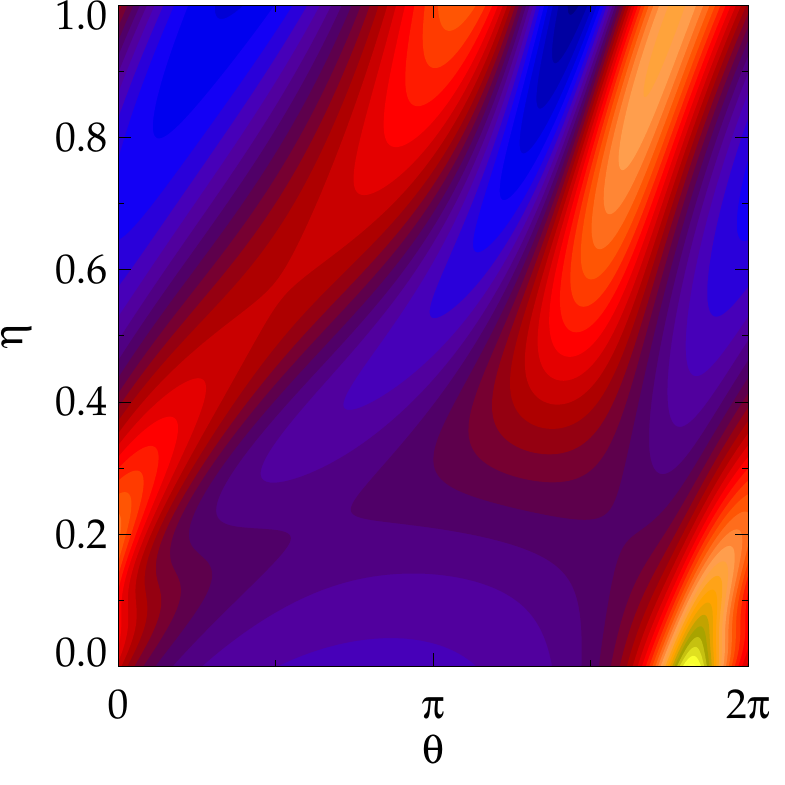}
  \includegraphics[width=0.40\textwidth]{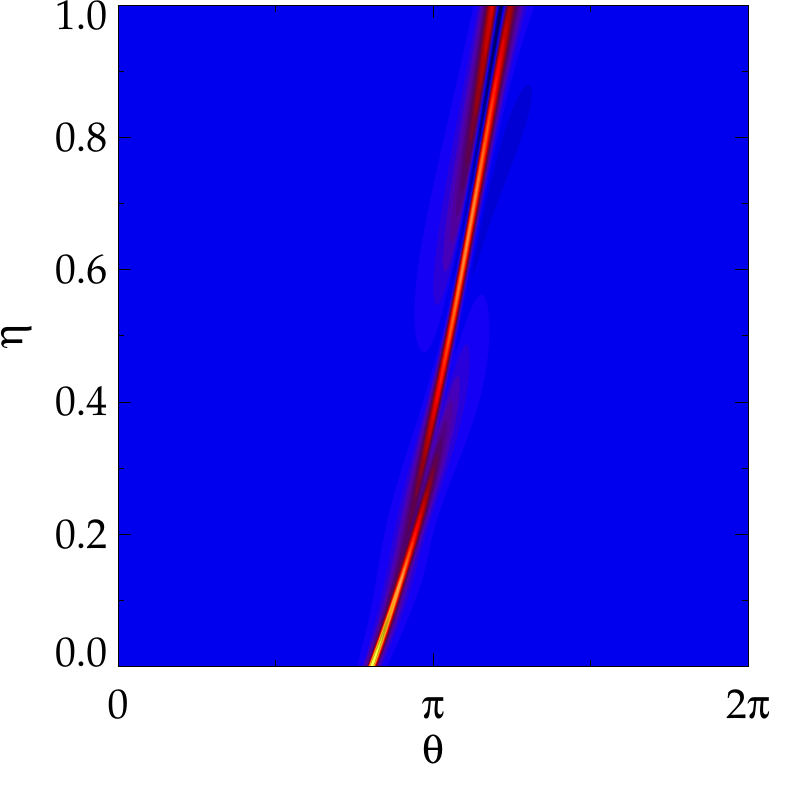}
  \caption{Trajectory $\mu_t(\theta,\mu)$ of \eqref{PDE} at time moments $t = 0, 3$ and 6
  (from top to bottom) computed for the optimal control input shown in Fig.~\ref{fig:ensemble_control}.
  The standard ``rainbow'' color table was used to code the isovalues: from black
  (minimal values), violet, \dots, to red (maximal values).}
  \label{fig:ensemble}
\end{figure}

\section{CONCLUSION}\label{sec-conc}

The goal of this paper is to present an approach based on the mean-field control paradigm
to solve problems of optimization and synchronization of oscillatory processes (here, the
addressed Theta model is among the simplest but prominent examples). The proposed
technique can be applied to any state-linear optimal control problem involving (finite or
infinite) non-interacting statistical ensembles of an arbitrary nature. In particular,
Algorithm~\ref{alg} can be easily adapted to some other neural model such as SNIPER model,
sinusoidal model etc. \cite{Li-TACON-2013}.   
 
We plan to continue this study in the way of natural generalization of model \eqref{Theta-k} by admitting the interaction between theta neurons,
\begin{equation*}
    \dot \theta_k = v_u(\theta_k, \eta_k)+\frac{1}{N}\sum_{j=1}^NK(\theta_k,\theta_j), \quad k=\overline{1,N},
\end{equation*}
where $K$ is certain interaction potential formalizing the spatial connectivity of neurons in the tissue.  This will result in control problems of the sort $(P_{1,2})$ stated over the nonlocal continuity equation
\[
\partial_t\mu_t + \partial_\theta\big(\left[v_{u} + K\star\mu_t\right]\, \mu_t\big) = 0
\]
involving the term
\[
(K\star\mu)(\theta) \doteq \int K(\theta,\zeta) \d \mu(\zeta).
\]
Such problems are not state-linear anymore, and the exact formula \eqref{incr} becomes inapplicable. For this case, a promising alternative could be an approach based on Pontryagin's maximum principle \cite{POGODAEV20203585}.

\bibliographystyle{IEEEtran}

\balance

\bibliography{Staritsyn}

\end{document}